\documentclass[12pt]{amsart}
\usepackage{amsmath,amssymb}

\newtheorem{thm}{Theorem}

\newtheorem{cor}[thm]{Corollary}

\newcommand{\R}{{\mathbb R}}

\newcommand{\modo}[1]{{\left|#1\right|}}
\newcommand{\normo}[1]{{\left\|#1\right\|}}

\newcommand{\snormo}[1]{{\mathopen\|#1\mathclose\|}}

\renewcommand{\div}{\text{\rm div}}

\begin{document}
\title[Finite time blow up for a Navier-Stokes like equation]
{Finite time blow up for a\\
Navier-Stokes like equation
}

\author{Stephen Montgomery-Smith}
\makeatletter
\address{Department of Mathematics\\
University of Missouri\\
Columbia, MO 65211}
\email{stephen@math.missouri.edu}
\urladdr{http://www.math.missouri.edu/\~{}stephen}
\thanks{The author was 
partially supported
by NSF grant DMS 9870026.}

\begin{abstract}
\noindent
We consider an equation similar to the Navier-Stokes equation.  
We show that
there is initial data that exists in every Triebel-Lizorkin or Besov space
(and hence in every Lebesgue and Sobolev space), such that after a finite
time, the solution
is in no Triebel-Lizorkin or Besov space (and hence 
in no Lebesgue or Sobolev space).  The purpose is to show the
limitations of the so called semigroup method for the Navier-Stokes
equation.
We also consider the possibility of existence of
solutions with initial data in the Besov space $\dot B^{-1,\infty}_\infty$.
We give initial data in this space for which there is no reasonable
solution for the Navier-Stokes like equation.
\end{abstract}

\maketitle

In this paper, we consider a simplified model for the 
Navier-Stokes equation --- what we call the cheap Navier-Stokes
equation.  For this equation, we show that for sufficiently
large initial data, that we get blow up in finite time.  
The purpose of this is {\em not\/} to indicate the 
possibility that the Navier-Stokes equation might blow up
in finite time --- indeed the author strongly believes
the opposite.  Rather, the purpose of this paper is to show limitations
in some of the methods used in studying the Navier-Stokes
equation.

Let us consider the following version of the Navier-Stokes equation:
\[
\frac{\partial u}{\partial t} = \Delta u - P(\div(u \otimes u)),
\]
where $t \mapsto u_t$ is an $\R^n$ valued function,
or tempered distribution, on $\R^n$.  
Here $P$ denotes the Leray projection that takes
a vector field to its divergence free part.
A tremendous amount of work has been done on
the very hard problem of determining if the solutions exist, if they
are unique, and to which spaces they belong.
One approach, the one we consider in this paper, 
is to consider mild solutions using what is often
called the semigroup approach.  
This is described in \cite{Ca1}, and
is used in many papers, for example, \cite{FK}, \cite{K}, \cite{GM}, \cite{KT}.
An example of this kind of result is that due to Kato \cite{K},
who showed that if $u_0 \in L_n(\R^n)^n$, 
then there is a solution $u \in C([0,T],(L_n(\R^n))^n)$
if either $T=T(u_0)$ is sufficiently 
small and $\snormo{u_0}_{L_n}$ is arbitrary or 
$\snormo{u_0}_{L_n}$ is small and $T=\infty$.
 
However, if one studies all these papers, one sees that they 
do not use all of the properties of the Navier-Stokes equation.  
Indeed, the methods seem to apply equally well to the the following
equation, the {\em cheap Navier-Stokes equation}:
\begin{equation*}
\frac{\partial u}{\partial t} = \Delta u + \sqrt{-\Delta} (u^2) .
\end{equation*}
Here $t\mapsto u_t$ is a scalar valued function or tempered distribution
on $\R^n$.
The semigroup approach is to consider
a mild solution of the cheap Navier-Stokes equation, that is, 
a solution to the equation 
$u = G(u)$,
where $G:C([0,T],X) \to C([0,T],X)$, 
with $X$ being a space of tempered distributions on $\R^n$:
\begin{equation*}
G(u)_t = 
e^{t \Delta} u_0 + \int_0^t e^{(t-s)\Delta} \sqrt{-\Delta}(u_s^2)) \, ds .
\end{equation*}
The method used 
to find a fixed point of $G$ is to show that $G$ is a contraction
mapping on $C([0,T],X)$ or on some subset of $C([0,T],X)$.
(Oftentimes one cannot obtain the result by working directly
with $C([0,T],X)$.  For example, in the case $X = L_n(\R^n)$ one
has to work with some space embedded into $X$, for as is shown by
Oru \cite{O}, the map $G$ is not even continuous on $C([0,T],L_n(\R^n))$.)

The semigroup
method by itself never seems to be able to obtain global results, that is,
results valid for both $T$ and the size of $u_0$ arbitrarily large.  
To obtain global results, one has to appeal to other
kinds of estimates for the Navier-Stokes equation, for example, energy
estimates (see, for example, \cite{KF} or \cite{Ca2}).

The purpose of this paper is to show that the semigroup method really is
limited in this way.  For if we could obtain
global results for the Navier-Stokes equation using only the semigroup
method, then the same methods would also
apply to the cheap Navier-Stokes equation.
This would then contradict the main result of this paper and its corollary.

Thus, if one is going to obtain global results for the Navier-Stokes
equation, one has to consider properties of the bilinear form
$P(\div(u \otimes u))$ that are not shared by the bilinear form
$\sqrt{-\Delta} (u^2)$.  Perhaps there is some mysterious cancellation
property in the first bilinear form that causes global regularity for
the Navier-Stokes equation.  In any case, the semigroup technique in
of itself, at least in the manner in which it has been applied to date,
is not going to solve the problem.

Let $e_1$ denote the first unit vector in $\R^n$, and for $x \in \R^n$,
$r>0$, let $B_r(x)$ denote $\{y:\modo{x-y}\le r\}$.
Let $\hat u$ denote the Fourier transform of $u$ with
respect to $x\in\R^n$.  

\begin{thm}
\label{blow-up}
Let $w$ be a tempered Schwartz function such that $\hat w$
is non-negative, has $L_1$ norm equal to $2$, and has support in 
$B_{1/4}(3e_1/4) \cup B_{1/4}(-3e_1/4)$.  (So $w$ is in every
Triebel-Lizorkin or Besov space.)
Then if $A > 2^{13/3}$, and if $u$ is a mild solution
to the cheap Navier-Stokes equation
whose Fourier transform is non-negative, with initial data $u_0 = A w$, 
then $u_t$
is not in any Triebel-Lizorkin or Besov space when $t = \log(2^{1/3})$.
\end{thm}

\begin{cor}
\label{no-Ln}
Let $w$ be as in Theorem~\ref{blow-up} (so $w \in L_n(\R^n)$), and let
$A > 2^{13/3}$.  Then there is 
no mild solution $u$ to the cheap 
Navier-Stokes equation, with $u_0=Aw$,
in $C([0,T],L_n(\R^n))$, for $T = \log(2^{1/3})$.
\end{cor}

We also present another result about the 
cheap Navier-Stokes
equation.  For the Navier-Stokes equation (and hence also for the
cheap Navier-Stokes equation), it turns out that the natural spaces
in which to consider solutions are of the form $C([0,T],X)$, where $X$
is {\em scale invariant}, that is,
$ \normo{\lambda f(\lambda \cdot)}_X = \normo f_X $
for all $\lambda \in (0,\infty)$.
Authors have considered which is the
largest scale invariant space for which one gets existence
results.  For example, recently Koch and Tataru \cite{KT} showed such
results when $X$ is the space of derivatives of functions with
bounded mean oscillation.

In fact, it can be shown (arguing similarly as in Frazier, Jawerth and Weiss 
in \cite{FJW} for the minimality of $B^{0,1}_1$)
that all scale invariant spaces of distributions,
that also contain all Schwartz functions, are
contained in the Besov space $\dot B^{-1,\infty}_\infty$.
Cannone \cite{Ca1} was able to obtain results for the Navier-Stokes
equation in the space $\dot B^{n/p-1,\infty}_p$, 
$n\leq p<\infty$, but left open the case $p=\infty$ corresponding
to the space $\dot B^{-1,\infty}_\infty$.

We will respond to this last case negatively for
the cheap Navier-Stokes equation.  Although this does not answer
the question for the Navier-Stokes equation, the author believes that
it is very possible that there is a similar non-existence result
for the Navier-Stokes equation.  But this would say more about the
nature of the space $\dot B^{-1,\infty}_\infty$ than about the 
Navier-Stokes equation itself.

\begin{thm}
\label{no-exist}
Let $w$ be a tempered Schwartz function whose Fourier transform
is non-negative, and has support in 
$B_{1}(0) \setminus B_{1/2}(0)$.
Let $v$ be the tempered distribution
\[ 
v(x) = \sum_{k=1}^\infty 2^k \cos((2^k-1) x_1) w(x) .
\]
Then $v \in \dot B^{-1,\infty}_\infty$, but there is
no mild solution $u\in C([0,T],{\mathcal S}')$ to the cheap 
Navier-Stokes equation, with non-negative
Fourier transform, and with $u_0=v$, for any $T>0$. (Here
${\mathcal S}'$ denotes the space of tempered distributions.)
\end{thm}

\noindent
The crucial observation for all these results is that
\begin{equation*}
\widehat{G(u)}_t(\xi) = 
e^{-t\modo{\xi}^2} \hat u_0 
+ \int_0^t e^{(s-t)\modo{\xi}^2} \modo{\xi} (\hat u_s*\hat u_s)(\xi) \, ds
\end{equation*}
(here $*$ denotes convolution).  Thus we see that if $\hat u, \hat v \ge 0$,
and that $\hat u \le \hat v$, then $\widehat{G(u)} \le \widehat{G(v)}$.

\begin{proof}[Proof of Theorem~\ref{blow-up}]
Using the usual embedding theorems, it may be seen that all of the
Triebel-Lizorkin and Besov spaces
embed into $\dot B_\infty^{a,\infty}$ for some $a \in \R$.
Let us recall the definition of the Besov space $\dot B_\infty^{a,\infty}$.
Let $\phi(x)$ be some 
tempered Schwartz function
whose Fourier transform is non-negative, 
whose support is ``mostly'' contained
in a band about $\modo\xi=1$, and such that
$\sum_{k=-\infty}^\infty \hat \phi_k$ is uniformly bounded above and below.
Here $\phi_k(x) = 2^{nk} \phi(2^k x)$.
Then $\dot B_\infty^{a,\infty}$ is the space of 
distributions on $\R^n$ for
which the norm
$
\snormo f_{\dot B_\infty^{a,\infty}} 
= \sup_k 2^{a k} \snormo{\phi_k* f}_{L_\infty} 
$
is finite.

Write $w_0$ for the function
$
\hat w_0 = \hat w I_{B_{1/4}(3e_1/4)} 
$,
and set $w_k = w_0^{2^k}$.  Since $w$ and $\hat w$ are real valued,
we see that $\hat w$ must be an even function.
Hence we quickly see that
$\hat w_k$ has $L_1$ norm equal to $1$, and is supported in 
$\{\xi:2^{k-1} \le \modo\xi \le 2^k\}$.

We will show by induction
that
$
\hat u_t \ge A^{2^k} \alpha_k(t) \hat w_k 
$
for $k>0$,
where $\alpha_k(t) = 2^{k-4(2^k-1)} e^{-2^k t} I_{t \ge t_k}$,
$t_0=0$, and $t_k = \log(2) \sum_{j=1}^k 2^{-2j}$.
The case $k=0$ follows since
$\hat u_t(\xi) \ge e^{-t \modo{\xi}^2} \hat u_0(\xi)$.
Suppose that our desired inequality holds for $k-1$.  Then
\begin{eqnarray*}
\hat u_t(\xi)
&\ge&
\int_0^t e^{(s-t)\modo{\xi}^2} \modo{\xi} (\hat u_s*\hat u_s)(\xi)) \, ds\\
&\ge&
\int_0^t e^{2^{2k}(s-t)} 2^{k-1} (A^{2^{k-1}} \alpha_{k-1}(s))^2
(\hat w_{k-1}*\hat w_{k-1})(\xi) \, ds\\
&\ge&
A^{2^k} 2^{k-4(2^k-1)} 2^{1+2k}
\int_{t_{k-1}}^t e^{2^{2k}(s-t)} e^{-2^k s} \, ds \, 
\hat w_k(\xi) \\
&\ge&
A^{2^k} 2^{k-4(2^k-1)} e^{-2^k t} 2^{1+2k} 
\int_{t_{k-1}}^t e^{2^{2k}(s-t)} \, ds \, \hat w_k(\xi)\\
&=&
A^{2^k} 2^{k-4(2^k-1)} e^{-2^k t}
2(1-e^{2^{2k}(t_{k-1}-t)}) \hat w_k(\xi)\\
&\ge&
A^{2^k} 2^{k-4(2^k-1)} e^{-2^k t} \hat w_k(\xi)
\end{eqnarray*}
whenever $1-e^{2^{2k}(t_{k-1}-t)} \ge 1/2$, that is, $t_k-t_{k-1} \ge 
2^{-2k} \log(2)$.

Let $t_\infty = \lim_{k\to\infty} t_k
= \log(2^{1/3})$.
Since 
$\hat \phi_k \hat u_{t_\infty} \ge A^{2^k} \alpha_k(t_\infty) \hat w_k \ge 0$,
\begin{eqnarray*}
\snormo{u_{t_\infty}}_{\dot B_{\infty}^{a,\infty}} 
&=& 
\sup_k 2^{ak} \snormo{\phi_k*u_{t_\infty}}_{L_\infty}
=
\sup_k 2^{ak} \snormo{\hat \phi_k \hat u_{t_\infty}}_{L_1} 
\\
&\ge&
\sup_k 
A^{2^k} 2^{ak+k-4(2^k-1)} e^{-2^k t_\infty} .
\end{eqnarray*}
It is clear that this is infinite if $A>2^4 e^{t_\infty} = 2^{13/3}$.
\end{proof}

\begin{proof}[Proof of Corollary~\ref{no-Ln}]
Suppose for a contradiction that there is a solution 
$u \in C([0,T],L_n(\R^n))$.  By the methods of \cite{K}, 
we know that there is a 
number $\epsilon>0$, depending only upon $\normo u$, such that
for every $t \in [0,T]$ that there is a mild solution 
$v \in C([t,t+\epsilon],L_n(\R^n))$ to the cheap Navier-Stokes equation, 
with $v_t = u_t$, that is obtained by starting with 
$v^{(0)}_s = e^{(s-t)\Delta} u_t$,
and iterating a function $G$ similar to that defined above.  
By the uniqueness result of Furioli, Lemari\'e-Rieusset 
and Terraneo \cite{FLT},
we have that $v_s = u_s$ for $t\le s \le t+\epsilon$.
From this it is clear that if $\hat u_t \ge 0$, then
$\hat u_s \ge 0$ for $t\le s \le t+\epsilon$.  Applying this argument
several times, we see that $\hat u_t \ge 0$
for $0 \le t \le T$.
Then by Theorem~\ref{blow-up}, $u_T$ is not in any
Triebel-Lizorkin space, and hence in particular is not in
$L_n(\R^n)$.
\end{proof}

\begin{proof}[Proof of Theorem~\ref{no-exist}]
Suppose we have a non-negative solution $u:[0,T] \to {\mathcal S}'$ to 
the cheap Navier-Stokes equation with $u_0 = v$.  It is clear that
\[
\hat u_t(\xi) \ge e^{-t\modo\xi^2} \hat v(\xi)
\ge \sum_{k=1}^\infty 2^k e^{-2^{2k} t} 
{\textstyle\frac12}
(\hat w(\xi+(2^k-1)e_1)+\hat w(\xi-(2^k-1)e_1)) .
\]  
Thus
\begin{eqnarray*}
\hat u_t(\xi)
&\ge&
\int_0^t e^{(s-t)\modo\xi^2} \modo\xi \hat (u_s*\hat u_s)(\xi) \, ds \\
&\ge&
\int_0^t \sum_{k=1}^\infty 
e^{4(s-t)} 2^{2k-1} e^{-2^{2k+1}s} (\hat w * \hat w)(\xi) \, ds \\
&=&
\left(\sum_{k=1}^\infty 2^{2k-1} e^{-4t} \int_0^t e^{(4-2^{2k+1})s} \, ds \right)
(\hat w * \hat w) (\xi)\\
&=&
\left(\sum_{k=1}^\infty \frac{2^{2k-1}}{2^{2k+1}-4} e^{-4t} 
(1-e^{(4-2^{2k+1})t})\right)
(\hat w * \hat w)(\xi) ,
\end{eqnarray*}
which is infinite on the support of $\hat w * \hat w$ for $t>0$, and hence
$u_t$ cannot be in ${\mathcal S}'$.  
\end{proof}

The author would like to thank Marco Cannone for many useful discussions
and helpful remarks.

\end{document}